\documentclass[12pt]{amsart}
\usepackage{amsmath,amsthm,amsfonts,amssymb}
\usepackage{amsmath, amsthm, amssymb, amscd, amsfonts}
\usepackage[all]{xy}
\usepackage{amssymb, amsthm, amsmath}
\usepackage[english]{babel}
\usepackage{amsfonts}

\newcommand{\ra}{\rightarrow}

\newcommand{\ot}{\otimes}
\newcommand{\mtc}{\mathcal}

\newcommand{\Lam}{\Lambda}

\newcommand{\al}{\alpha}
\newcommand{\eps}{\epsilon}
\newcommand{\bn}{\begin}

\newcommand{\sub}{\subsection}

\newcommand{\D}{\Delta}

\numberwithin{equation}{section}

\newcommand{\dw}{\downarrow}
\newcommand{\uw}{\uparrow}

\newcommand{\ch}{\chi}
\newcommand{\mtr}{\mathrm}

\numberwithin{equation}{section}
\title[Clifford theory]
{Clifford theory for cocentral extensions}
\author{Sebastian  Burciu}
\address{Inst.\ of Math.\ ``Simion Stoilow" of the Romanian Academy
P.O. Box 1-764, RO-014700, Bucharest, Romania }
\email{sebastian.burciu@imar.ro}
\date{December 12, 2008}

\begin{document}
\begin{abstract} The classical Clifford correspondence for normal subgroups is considered in the more general setting of semisimple Hopf algebras. We prove that this correspondence still holds if the extension determined by the normal Hopf subalgebra is cocentral.
\end{abstract}
\maketitle
\section*{Introduction}
The starting point for Clifford theory is Clifford's paper \cite{C1} on representations of normal groups. Since then a lot of literature was written on the subject. Parallel theories for graded rings and Lie algebras were developed in \cite{dade} and \cite{Li1} respectively, as well as in other papers. A unifying setting for these theories was developed by Schneider \cite{Schgal} for Hopf Galois extensions. The main problem with this more general theory is that usually the stabilizer is not a Hopf subalgebra and is not an extension of the based ring.

A more general approach  was considered by Witherspoon in \cite{Sc1} for any normal extension of semisimple algebras. With a certain definition of the stabilizer it was proven in \cite{Sc1} that the Clifford correspondence holds.

In this paper we address an analogue of initial's Clifford approach for groups. We consider an extension of Hopf algebras $A/B$ where $B$ is a normal Hopf subalgebra of $A$ and let $M$ be an irreducible $B$-module. The conjugate $B$-modules of $M$ are defined as in \cite{coset} and the stabilizer $Z$ of $M$ is a Hopf subalgebra of $A$ containing $B$. We say that the Clifford correspondence holds for $M$ if induction from $Z$ to $A$ provides a bijection between the sets of $Z$ (respectively $A$)-modules that contain $M$ as a $B$-submodule.

Since $B$ is normal in $A$ also in the sense of \cite{Sc1} the results from this paper can also be applied.  It is shown that the Clifford correspondence holds for $M$ if and only if $Z$ is a stabilizer in the sense proposed in \cite{Sc1}. A necessary and sufficient condition for this to happen is given in Proposition \ref{num}.  Our approach uses the character theory for Hopf algebras and normal Hopf subalgebras. If the extension $$\begin{CD}k @ > >> B  @> i >>  A @ > \pi >> H@ > >> k \end{CD}$$ is cocentral then we prove that this condition is satisfied (see Corollary \ref{coc}). Recall that a such extension is cocentral if $H^* \subset \mtc{Z}(A^*)$ via $\pi^*$.

The paper is organized as follows. First section recalls the character theory results for Hopf subalgebras that are further needed. The next section defines the conjugate module and introduces the stabilizer as a Hopf subalgebra. The necessary and sufficient condition for the Clifford correspondence to hold is proven in this section. Third section considers the case when the quotient Hopf algebra is a finite group algebra. A different approach gives in these settings another criterion for the Clifford correspondence to hold (see Theorem \ref{princ}). As a corollary of this it is proven that the Clifford correspondence holds for cocentral extensions. In the last section of the paper a counterexample of a non cocentral extension where the Clifford correspondence does not hold anymore is given.

For a vector space $V$ over of a field $k$  by $|V|$ is denoted the dimension $\mtr{dim}_kV$. The comultiplication, counit and antipode of a Hopf algebra are  denoted by $\Delta$, $\epsilon$ and
$S$, respectively. We use Sweedler's notation $\D(x)= x_1\ot x_2$ for all $x\in H$ with the sum symbol dropped.  All the other notations are those used in \cite{Montg}.
All considered modules are left modules.


\section{Normal Hopf subalgebras}

Throughout of this paper $A$ will be a finite dimensional semisimple Hopf algebra over an algebraically closed field $k$ of characteristic zero. Then $A$ is also cosemisimple and $S^2=\mtr{Id}$ \cite{Lard}. The set of irreducible characters of $A$ is denoted by $\mtr{Irr}(A)$. The Grothendieck group $\mtc{G}(A)$ of the category of finite dimensional left $A$-modules is a ring under the tensor product of modules. Then $C(A)=\mtc{G}(A)\ot_{\mathbb{Z}}k$ is a semisimple subalgebra of $A^*$ \cite{Zhc} and it has a basis given by the characters of the irreducible $A$-modules.

Let $B$ be a Hopf subalgebra $A$. By Corollary 2.5 of \cite{coset} there is a coset decomposition for $A$
$$A=\oplus_{C/\sim}BC.$$ where $\sim$ is an equivalence relation on the set of simple subcoalgebras of $A$ given by $C \sim C'$ if and only if $BC=BC'$. In \cite{coset} this equivalence relation is denoted by $r^{A}_{ _{B,\;k}}$.

Since $A$ is also cosemisimple \cite{Lard} the set of simple subcoalgebras of $A$ is in bijection with the set of irreducible characters of $A^*$ (see \cite{Lar} for this correspondence).

Suppose now that $B$ is a normal Hopf subalgebra of $A$. Recall that this means $a_1BS(a_2)\subset B$ for all $a \in A$. If $\ch $ and $\mu$ are two irreducible characters of $A$ it can be proven that their restriction to $B$ either have the same irreducible constituents or they don't have common constituents at all. Define  $\ch \sim \mu$ if and only $m_{ _B}(\ch \dw_B^A,\;\mu\dw_B^A)>0$. With the above notations this is the  equivalence relation $r^{A^*}_{ _{H^*,\;k}}$ for the inclusion $H^* \subset A^*$ where $H=A//B$ is the quotient Hopf algebra.

Let $\mtc{A}_1,\;\cdots,\mtc{A}_l$ the equivalence classes of the above relation and $$a_i=\sum_{\ch \in \mtc{A}_i}\ch(1)\ch$$ for $1 \leq i \leq l$.

This equivalence relation determines an equivalence relation on the set of irreducible characters of $B$. Two irreducible $B$-characters $\al $ and $\beta$ are equivalent if and only if they are constituents of  $\ch\dw_B^A$ for some irreducible character $\ch$ of $A$.

Let $\mtc{B}_1,\cdots,\mtc{B}_l$ be the equivalence classes of this new equivalence relation and let $$b_i=\sum_{\al \in \mtc{B}_i}\al(1)\al.$$ The induction-restriction formulae from \cite{coset} can be written
as
\bn{equation}\label{rch}
\frac{\ch\dw_B^A}{\ch(1)}=\frac{b_i}{b_i(1)}
\end{equation} and
\bn{equation}\label{ind}
\frac{\al\uw_{ _B}^A}{\al(1)}=\frac{|A|}{|B|}\frac{a_i}{a_i(1)}
\end{equation} if $\ch \in \mtc{A}_i$ and $\al \in \mtc{B}_i$.

Since the regular character of $A$ restricts to $\frac{|A|}{|B|}$ copies of the regular character of $B$ it follows that

\bn{equation}
\label{rall}a_i\dw_B^A=\frac{|A|}{|B|}b_i.
\end{equation}

In particular $a_i(1)=\frac{|A|}{|B|}b_i(1)$ for all $1 \leq i \leq l$(see also 4.1 of \cite{coset}).
\section{Conjugate modules and stabilizers}
Let  $M$ be an irreducible $B$-module with character $\al \in C(B)$. We recall the following notion of conjugate module introduced in \cite{coset}. It was also previously considered in \cite{Schgal} in the cocommutative case.

If $W$ is an $A^*$-module then $W\ot M$ becomes a $B$-module with

\bn{equation}\label{def}
 b(w\ot m)=w_0 \ot (S(w_1)bw_2)m
\end{equation}

Here we used that any left $A^*$-module $W$ is a right
$A$-comodule via $\rho(w)=w_0\ot w_1$.
It can be checked that if $W \cong W'$ as $A^*$-modules then $W\ot M \cong W'\ot M$. Thus for any irreducible character $d \in \mtr{Irr}(A^*)$ associated to a simple $A$-comodule $W$ one can define the $B$-module $\;^dM \cong W\ot M$. If $\al$ is the character of $M$ then the character $^d\al$ of $^dM$ is given by

\bn{equation}\label{chfom}
^d\al(x)=\al(Sd_1xd_2)
\end{equation}
for all $x \in B$ (see Proposition 5.3 of \cite{coset}).

\bn{rem}
\label{eqconst}
From Proposition 5.12 of \cite{coset} it follows that the equivalence class of a character $\al \in \mtr {Irr}(B)$ is given by all the irreducible constituents of  $^d\al$ as $d$ runs through all irreducible characters of $H^*$.
\end{rem}

Fix $\al \in \mtr{Irr}(B)$ and suppose that $\al \in \mtc{B}_{i}$ for some index $i$.

\bn{prop}
The set $\{d \in \mtr{Irr}(A^*)\;| \;^d\al=\eps(d)\al\}$ is closed under multiplication and $``\;^*\;"$. Thus it generates a Hopf subalgebra $Z$ of $A$ that contains $B$.
\end{prop}

\bn{proof}
Since $^d(^{d'}\al)=\; ^{dd'}\al$ it follows that the above set is closed under multiplication. Since  $d^*$ is a constituent of some power of $d$  it also follows that the set is closed under $``\;^*\;"$ too. Thus it generates a Hopf subalgebra $Z$ of $A$ (see \cite{NR}) with $Z=\oplus_{C} C$ where the sum is over all simple subcoalgebras of $H$ whose irreducible characters $d$ satisfy $^d\al=\eps(d)\al$. If $d \in B$ then $^d\al(x)=\al(Sd_1xd_2)=\al(xd_2S(d_1))=\eps(d)\al(x)$ for all $x\in B$. Therefore $B \subset Z$.
\end{proof}

$Z$ will be called the stabilizer of $\al$ in $A$.
\bn{rem}\label{tpr}
If $C$ is any subcoalgebra of $H$ then $C\ot M$ has a structure of $B$-module as above using the fact that $C$ is a right $A$-comodule via $\D$. Then $C\ot M \cong M^{|C|}$ as $B$-modules if and only if $C \subset Z$. \end{rem}
\bn{rem}
If $A=kG$ and $B=kN$ for a normal subgroup $N$ then $Z$ coincides with the stabilizer of $\al$ introduced in \cite{C1}.
\end{rem}
\subsection{On the stabilizer}
Since $B$ is normal in $Z$ one can define as above two equivalences relations, on
$\mtr{Irr}(Z)$ respectively  $\mtr{Irr}(B)$. Let ${\mtc{Z}_1}, \cdots ,\mtc{Z}_r$ be the equivalence classes in $\mtr{Irr}(Z)$ and  $\mtc{B'}_1, \cdots ,\mtc{B'}_r$ be the corresponding equivalence classes in $\mtr{Irr}(B)$. 

Remark \ref{eqconst} implies that $\al$ by itself form an equivalence class of $\mtr{Irr}(B)$, say $\mtc{B'}_1$. Then clearly the corresponding equivalence class ${\mtc{Z}_1}$ is given by $${\mtc{Z}_1}=\{\psi \in \mtr{Irr}(Z)|\;\psi\dw_{ _B}
\;\textnormal{contains}\; \al\}.$$
 Formula \ref{rch} becomes in this situation $\psi\dw_{ _B}^Z=\frac{\psi(1)}{\al(1)}\al$ for all
$\psi \in {\mtc{Z}_1}$.  Let $\psi_{\al}=\sum_{\psi \in  {\mtc{Z}_1}}\psi(1)\psi$.
Then $\psi_{\al}\dw_{ _B}^Z=\frac{|Z|}{|B|}\al(1)\al$ by \ref{rall} and
$\psi_{\al}(1)=\frac{|Z|}{|B|}\al(1)^2$.

\bn{lemma}\label{infr}
With the above notations $${\psi}_{\al}\uw^A_Z=\frac{\al(1)^2}{b_i(1)}a_i.$$
\end{lemma}

\bn{proof}
One has  $\al\uw_B^Z=\frac{|Z|}{|B|}\frac{\al(1)}{{\psi}_{\al}(1)}{\psi}_{\al}$ by \ref{ind}. But ${\psi}_{\al}(1)=\frac{|Z|}{|B|}\al(1)^2$ and the last formula becomes $$\al\uw_B^Z=\frac{{\psi}_{\al}}{\al(1)}$$ Thus $\al\uw_B^A=(\al\uw_B^Z)\uw_Z^A=\frac{{\psi}_{\al}\uw^A_Z}{\al(1)}$. On the other hand $\al\uw_B^A=\frac{|A|\al(1)}{|B|a_i(1)}a_i$ and one gets that:

$$ {\psi}_{\al}\uw^A_Z=\frac{|A|\al(1)^2}{|B|a_i(1)}a_i=\frac{\al(1)^2}{b_i(1)}a_i.$$
\end{proof}

\subsection{ Definition of the Clifford correspondence}The above Lemma implies that for any $\psi \in \mtc{Z}_1$  all the irreducible constituents of  $\psi\uw^A_Z$ are in $\mtc{A}_{i}$. We say that the Clifford correspondence holds for the irreducible character $\al \in \mtc{B}_i$ if $\psi\uw^A_Z$ is irreducible for any irreducible character $\psi \in \mtc{Z}_1$ and the induction function $$\mtr{ind}: {\mtc{Z}_1}\ra \mtc{A}_{i}$$ given by $\mtr{ind}(\psi)=\psi\uw^A_Z$ is a bijection.
\subsection{Clifford theory for normal subrings \cite{Sc1}}

Let $B \subset A$ an extension of $k$-algebras. An ideal $J$ of $B$ is called $A$-invariant if $AJ=JA$. Following \cite{Sc1} the extension $A/B$ is called normal if every two sided ideal of $B$ is $A$-invariant. Witherspoon gave a general Clifford correspondence for normal extensions. Let $M$ be a $B$-module. Then $M$ is called $A$ stable if the module $M\uw_B^A \dw_B^A $ is isomorphic to direct sum of copies of $M$.
A stabilizer $S$ of $M$ is a semisimple algebra $S$ such that $B \subset S \subset A$, $B$ is a normal subring of $S$, $M$ is $S$-stable, and $M -\mtr{soc} (M\uw_B^A\dw_B^A)=M -\mtr{soc} ( M\uw_B^S\dw_B^S)$. Here the $M$-socle of a $B$ -module is the sum of all its submodules isomorphic to $M$.

Next we investigate a relationship between the stabilizer $Z$ previously defined and the notion of stabilizer defined as above for normal extensions. It is easy to see that if $B$ is a normal Hopf subalgebra of $A$ then the extension $A/B$ is normal in the above sense (see also Proposition 5.3 of \cite{Sc1}). By the same argument $B$ is normal in $Z$ and from Remark \ref{eqconst} it follows that $M$ is $Z$-stable. Thus $Z$ is a stabilizer in the above sense if and only if the socle condition is satisfied. In terms of characters this can be written as  $m_{ _B}(\al\uw_B^Z\dw_B^Z,\;\al) = m_{ _B}(\al\uw_B^A\dw_B^A,\;\al)$ where $\al$ is the character of $M$.
\bn{prop}\label{num}With the above notations:
\bn{enumerate}
\item  $|Z| \leq \frac{|A|\al(1)^2}{b_i(1)}$.

\item  Equality holds if and only if $Z$ is a stabilizer in the sense of \cite{Sc1}.
\end{enumerate}
\end{prop}

\bn{proof}
Clearly
$m_{ _B}(\al\uw_B^Z\dw_B^Z,\;\al) \leq m_{ _B}(\al\uw_B^A\dw_B^A,\;\al)$ and  equality holds if and only if $Z$ is a stabilizer in the sense of \cite{Sc1}.

Let as before $s=\frac{|A|}{|B|}$ be the index of $B$ in $A$ and $s'=\frac{|Z|}{|B|}$ be the index of $B$ in $Z$.

Using  formulae \ref{ind} and \ref{rall} it can be seen that :
\bn{eqnarray*}
m_{ _B}(\al\uw_B^A\dw_B^A,\;\al)& \!=\! &  \frac{s\al(1)}{a_i(1)}m_{ _B}(a_i\dw^A_B,\;\al)
=\frac{s^2\al(1)}{a_i(1)}m_{ _B}(b_i,\;\al)
\\ & = & \frac{\al(1)^2s^2}{a_i(1)} \\ & = & \frac{\al(1)^2}{b_i(1)}s
\end{eqnarray*}
A similar argument applied to the extension $B \subset Z$ gives
\bn{equation*}
m_{ _B}(\al\uw_B^Z\dw_B^Z,\;\al)=\frac{\al(1)^2s'}{b'_1(1)}=s'
\end{equation*}
since in this situation $b'_1=\al(1)\al$. Thus $s'=\frac{|Z|}{|B|}\leq  s\frac{\al(1)^2}{b_i(1)}=\frac{|A|}{|B|}\frac{\al(1)^2}{b_i(1)}$ which gives the required inequality.
\end{proof}

\bn{rem}
If $A=kG$ and $B=kN$ for a normal subgroup $N$ then the above inequality is equality. It states that the number of conjugate modules of $\al$ is the index of the stabilizer of $\al$ in $G$.
\end{rem}
\subsection{Clifford correspondence}
\bn{thm}\label{main}
The Clifford correspondence holds for $\al$ if and only if  $Z$ is a stabilizer in the sense given in \cite{Sc1}.
\end{thm}

\bn{proof}
If $Z$ is a stabilizer in the sense given in \cite{Sc1} then the Clifford correspondence holds by Theorem 4.6 of the same paper.

Conversely, suppose that the map $$\mtr{ind}: {\mtc{Z}_1}\ra \mtc{A}_i$$ given by $\mtr{ind}(\psi)=\psi\uw^A_Z$ is a bijection. Thus for any $\psi \in \mtc{Z}_1$ there is a $\ch \in \mtc{A}_i$ such that $\psi\uw^A_Z=\ch$. Note that this implies $\psi(1)=\frac{|Z|}{|A|}\chi(1)$.

Since $\mtr{ind}$ is a bijection one can write
$${{\psi}_{\al}}\uw^A_Z=\sum_{\psi \in {\mtc{Z}_1}}\psi(1)\psi\uw^A_Z=\sum_{\ch \in \mtc{A}_i}\frac{|Z|}{|A|}\chi(1)\ch=\frac{|Z|}{|A|}a_i$$
which implies that ${\psi}_{\al}(1)=(\frac{|Z|}{|A|})^2a_i(1)$. Lemma \ref{infr} implies ${\psi}_{\al}(1)=\frac{|Z|}{|A|}\frac{\al(1)^2}{b_i(1)}a_i(1)$ and therefore one gets  $|Z|= \frac{|A|\al(1)^2}{b_i(1)}$. Proposition \ref{num} implies that $Z$ is a stabilizer in the sense given in \cite{Sc1}.
\end{proof}

\section{Extensions of Hopf algebras}
Let $B$ be a normal Hopf subalgebra of $A$ and $H=A//B$. Then we have the extension
\bn{equation}\begin{CD}\label{ext}k @ > >> B  @> i >>  A @ > \pi >> H@ > >> k \end{CD}\end{equation}
and $A/B$ is an $H$-Galois extension with the comodule structure $\rho:A \ra A\ot H$ given by $\rho=(\mtr{id} \ot \pi)\D$.

\bn{rem}\label{restric}
The restriction functor form $A$-modules to $B$-modules induces a map $\mtr{res}:C(A)\ra C(B)$. It is easy to see that $\mtr{res}=i^*|_{C(A)}$, the restriction of $i^*:A^* \ra B^*$  to the subalgebra of characters $C(A)$. By duality, $\pi|_{C(A^*)}$ is the restriction map of $A^*$-characters to $H^*$ (here $H^* \subset A^*$ via $\pi^*$).
\end{rem}

\subsection{Results on Hopf Galois extensions}
In this subsection we recall few facts about Clifford theory for Hopf Galois extensions over finite group algebras $H=kF$ from \cite{Schgal}. (see also \cite{dade}). Let $A/B$ be a Hopf Galois extension over $H=kF$ via the comodule map $\rho: A \ra A\ot kF$. For any $f \in F$ let $A_f=\rho^{-1}(A \ot kf)$. Since $A/B$ is a Hopf Galois extension one has that $A=\oplus_{f \in F}A_f$ is a strongly graded algebra by $F$ with  $A_1=B$. The functor $A_f \ot_B -: _B\mtc{M}\ra _B\mtc{M}$ is an equivalence of categories since $A_f$ is an invertible $B$-bimodule, $A_f \ot_B A_{f^{-1}}=B$.
In particular for any simple $B$-module $M$ then $A_f \ot_BM$ is also a simple $B$-module. From this it follows that the group $F$ acts on the irreducible representations of $B$ by $f.M:=A_f\ot_BM$

The stabilizer $H$ of $M$ is defined as the set of all $f \in F$ such that $A_f \ot_BM\cong M$ as $B$-modules. It is a subgroup of $F$. Let $S:=A(H)=\rho^{-1}(A \ot kH)=\oplus_{h \in H}A_h$. Then the induction map $\mtr{ind}$  $$ \{V \in S-mod:  V\dw_B^S contains\; M \}\ra\\ \{P \in A-{mod}: P\dw_B^A contains  \;M \}$$ given by $\mtr{ind}(M)=S\ot_BM$ is a bijection.

\subsection{Extensions by $kF$}
For the rest of this section we suppose that $H=kF$ for some finite group $F$. Then $H^*=k^F$ is a normal Hopf subalgebra of $A^*$ and one can define the same equivalence relations form the beginning of this paper for this extension. Since $\mtr{Irr}(k^F)=F$ this gives a partition of the group $F=\bigsqcup_{j=1}^m\mtc{F}_j$. Then by Remark \ref{restric} formula \ref{rall} applied to this situation implies that for any $d \in \mtr{Irr}(A^*)$ there is an unique index $j$ such that
\bn{equation}\label{rF}\pi(d)=\frac{\eps(d)}{|\mtc{F}_j|}\sum_{f \in \mtc{F}_j}f.\end{equation}

\subsection{Dimension of the orbit}
Let $M$ be an irreducible representation of $B$ with character $\al$ and let $H \leq F$ be the stabilizer of $M$. Since $|f.M|=|M|$ it follows that all the irreducible representations in the equivalence class of $M$ have the same dimension. If $s$ is their number then clearly $s=\frac{|F|}{|H|}$. Suppose now that $\mtc{B}_i$ is the equivalence class of $\al$. The above results implies that  $\mtc{B}_i$ coincide with the set of characters of the irreducible modules $f.M$. Thus
\bn{equation}\label{orbit}
b_i(1)=s\al(1)^2=\frac{|F|}{|H|}\al(1)^2.
\end{equation}

\subsection{Coset decomposition}

Recall the coset decomposition for $A$
\bn{equation}\label{cstform}
A=\oplus_{C/\sim}BC.
 \end{equation}
where $\sim$ is an equivalence relation on the set of simple subcoalgebras of $A$ given by $C \sim C'$ if and only if $BC=BC'$. Note that $BC=CB$ for any simple subcoalgebra $C$ of $A$ since $B$ is a normal Hopf subalgebra (see also \cite{coset}).

\bn{lemma}\label{image}
Suppose that $\pi:A \ra kF$ is a surjective map of Hopf algebras where $F$ is a finite group. Let $C$ be a simple subcoalgebra of $A$ with irreducible character $d$ and suppose $\pi(d)=\sum_{g \in \mtc{A}}a_gg$ where $\mtc{A}\subset F$ and $a_g$ are positive integers for all $g \in \mtc{A}$. Then $\pi(C)=\oplus_{g \in \mtc{A}}kg$.
\end{lemma}

\bn{proof}
$\pi(C)$ is a subcoalgebra of $kF$ and therefore $\pi(C)=\oplus_{g \in \mtc{B}}kg$. It is enough to show $\mtc{A}=\mtc{B}$. Clearly $\mtc{A}\subset \mtc{B}$. For any $g \in \mtc{B}$ let $k_g$ be a copy of the field $k$. By duality $\pi^*$ induces an embedding of the semisimple algebra $R=\prod_{g \in \mtc{B}}k_g$ in the matrix algebra $C^*=M_{\eps(d)}(k)$. Writing the primitive idempotents of $R$ in terms of the primitive idempotents of $C^*$ it follows that $\mtc{B}\subset \mtc{A}$.
\end{proof}

\bn{lemma}\label{cdec}
Assume that $H=kF$ for some finite group $F$. Let $d \in \mtr{Irr}(A^*)$ associated to the simple subcoalgebra $C$. If $$\pi(d)=\frac{\eps(d)}{|\mtc{F}_j|}\sum_{f \in \mtc{F}_j}f$$ then the coset $BC= \oplus_{f \in \mtc{F}_j}A(f)$.
\end{lemma}

\bn{proof}
Let $A_s=\oplus_{f \in \mtc{F}_s}A(f)$ for all $1 \leq s \leq m$. Then $A=\oplus_{s=1}^mA_s$.
The above lemma implies that $\pi(C)=\oplus_{f \in \mtc{F}_j}kf$. Since $\pi(BC)=\pi(C)$ this shows $BC \subset  A_j$. The coset decomposition formula \ref{cstform} forces $BC =  A_j$.
\end{proof}

\bn{thm}\label{princ}
Suppose that $H=kF$ for some finite group $F$. Let $M$ be an irreducible representation of $B$ with character $\al$ and let $H \leq F$ be the stabilizer of $M$. Then $Z\subset S:=A(H)$ and the Clifford correspondence holds for $\al$ if and only if $Z=S$.
\end{thm}

\bn{proof}
Since $A/B$  is an Hopf Galois extension over $H= kF$ it follows as above that $A$ is strongly $F$-graded with $A=\oplus_{f \in F}A_f$. First we will show that $Z\subset S=A(H)$. Recall the definition of $Z$ as the sum of all simple subcoalgebras $C$ whose irreducible characters $d$ verify the property $^d\al=\eps(d)\al$. Let $C$ be such an algebra with character $d$. As above there is a $j$ such that $\pi(d)=\frac{\eps(d)}{|\mtc{F}_j|}\sum_{f \in \mtc{F}_j}f$. It is easy to see that the canonical map $CB \ot M \ra CB\ot_BM $ is a morphism of $B$-modules. Since $BC \subset Z$ is a subcoalgebra Remark \ref{tpr} implies $CB \ot M \cong M^{|CB|}$ as $B$-modules. Thus $CB\ot_BM$ is a sum of copies of $M$. By Lemma \ref{cdec} $BC\ot_BM=\oplus_{f \in \mtc{F}_j}A(f)\ot_BM$ which shows that $\mtc{F}_j\subset H$ and therefore $C \subset A(H)$ by Lemma \ref{image}. Thus $Z=\sum_{C \subset Z}C \subset A(H)$.

Since $S/B$ is a $kH$-Hopf Galois extension it follows that $|S|=|B||H|$. Using formula \ref{orbit} it follows that $|S| = \frac{|A|\al(1)^2}{b_i(1)}$ if $\al \in \mtc{B}_i$. Then Theorem \ref{main} shows that the Clifford correspondence holds if and only if $|Z| = |S|$.
\end{proof}

It is easy to see that $\D_A(S)\subset A\ot S$.

\bn{cor}
Suppose that $H=kF$ for some finite group $F$. Let $M$ be an irreducible representation of $B$ with character $\al$ and let $H \leq F$ be the stabilizer of $M$. Then the Clifford correspondence holds for $\al$ if and only $S$ is a Hopf subalgebra of $A$.
\end{cor}

\bn{proof}
Any $S$-module which restricted to $B$ contains $M$ is a direct sum of copies of $M$ as a $B$-module by Corollary 2.2 of \cite{Schgal}. If $S$ is a Hopf algebra then Remark \ref{eqconst} applied to the extension $S/B$ implies that $S \subset Z$. Thus $S=Z$.
\end{proof}

\bn{cor}\label{coc}
Suppose that the extension \ref{ext} is cocentral. Then the Clifford correspondence holds for any irreducible $B$-module $M$.
\end{cor}

\bn{proof}
Since $H^*$ is commutative there is a finite group $F$ such that $H=kF$.
It is easy to see that $H^* \subset \mtc{Z}(A^*)$  via $\pi^*$ if and only if $\pi(a_1) \ot a_2=\pi(a_2)\ot a_1$ for all $a \in A$. This last relation implies that $S$ is a Hopf subalgebra of $A$ and the previous corollary finishes the proof.
\end{proof}

\section{A Counterexample}

Let $\Sigma = FG$ be an an exact factorization of finite groups. This gives a right action  $\lhd : G \times F \ra G$ of $F$ on the set $G$, and a left action $\rhd : G \times F \ra F$ of $G$ on the set $F$ subject to the following two conditions:
$$ s\rhd xy = (s \rhd x)((s\lhd x) \rhd y)\;\;\;\;\;
 st\lhd x = (s\lhd (t \rhd x))(t \lhd x)$$

The actions $\rhd$ and $\lhd$ are determined by the relations $gx = (g\rhd x)(g \lhd x)$ for all
$ x \in F$, $g \in G$. Note that $1 \rhd x=x$ and $s\lhd1=s$.

Consider the Hopf algebra $A=k^G \# kF$ \cite{Masext} which is a smashed product and coproduct using the above two action. The structure of $A$ is given by:

$$(\delta_gx)(\delta_hy) = \delta_{g\lhd x,h} \delta_gxy$$
$$\D(\delta_gx) = \sum_{
st=g} \delta_s(t \rhd x) \otimes \delta_tx$$

Then $A$ fits into the abelian extension
\bn{equation}\begin{CD}\label{abext}k @ > >> k^G  @> i >>  A @ > \pi >> kF@ > >> k \end{CD}\end{equation}

As above $F$ acts on $\mtr{Irr}(k^G)=G$. It is easy to see that this action is exactly $\lhd$.
Let $g \in G$ and $H$ be the stabilizer of $g$ under $\lhd$. Using the above notations it follows that $S=A(H)=k^G \# kH$. We will construct an example where $S$ is not a Hopf algebra and therefore the Clifford correspondence does not hold $g \in G$. Remark that the above comultiplication formula implies $S$ is a Hopf subalgebra if and only if $G \rhd H \subset H$.

Consider the exact fact factorization $\mathbb{S}_4=\mathbb{C}_4\mathbb{S}_3$ where $\mathbb{C}_4$ is generated by the four cycle $g=(1234)$ and $\mathbb{S}_3$ is given by the permutations that leave $4$ fixed. If $t=(12)$ and $s=(123)$ then the actions $\lhd$ and $ \rhd$ are given in Tables 1 and 2.

\begin{table}[t]
\begin{center}
\tiny{\begin{tabular}{|p{1,5cm}|p{1cm}|p{1cm}|p{1cm}|}

\hline {$\mathbb{C}_4 \lhd \mathbb{S}_3$ } &  {\bf $g$\newline} & {\bf $g^2$
} & {\bf $g^3$}
\\ \hline $t\newline$ & $g$  & $g^3$    & $g^2$
\\ \hline  $s\newline$ & $g^2$ & $g^3$  &$g$
 \\ \hline  $s^2\newline$ & $g^3$  & $g$ & $g^2$
 \\ \hline  $st\newline$ & $g^3$ &$g^2$  &$g$
\\ \hline $ts\newline$  & $g^2$  &  $g$ & $g^3$
\\ \hline
\end{tabular}}
\end{center}
\caption{The right action of $\mathbb{S}_3$ on $\mathbb{C}_4$}\label{tabl}
\end{table}

\begin{table}[t]
\begin{center}
\tiny{\begin{tabular}{|p{1,5cm}|p{1cm}|p{1cm}|p{1cm}|p{1cm}|p{1cm}|}

\hline {$\mathbb{C}_4\rhd \mathbb{S}_3$} &  {\bf $t$\newline} & {\bf $s$
} & {\bf $s^2$} & {\bf $st$\newline} & {\bf $ts$\newline}
\\ \hline $g\newline$ & $ts$ &  $t$  & $s$  & $st$ & $s^2$
\\ \hline  $g^2\newline$ & $s^2$ & $ts$   & $t$ & $st$ & $s$
 \\ \hline  $g^3\newline$ & $s$ & $s^2$ & $ts$   & $st$ & $t$

\\ \hline
\end{tabular}}
\end{center}
\caption{The left action of $\mathbb{C}_4$ on $\mathbb{S}_3$}\label{tabl}
\end{table}
 The stabilizer of the element $g$ is the subgroup $\{1, t\}$ which is not invariant by the action of $\mathbb{C}_4$. Thus the Clifford correspondence does not hold for $g$.
\bibliographystyle{amsplain}
\bibliography{cocentral}
\end{document}